\documentclass[11pt]{amsart}
\usepackage{amsthm}
\usepackage{amsmath}
\usepackage{amssymb}
\usepackage{amsfonts}
\usepackage{amscd}
\usepackage[mathscr]{eucal}
\usepackage{verbatim}
\usepackage{latexsym}
\usepackage{graphics}
\usepackage[all]{xy}

\theoremstyle{plain}
\newtheorem{thm}{Theorem}[section]
\newtheorem{lem}{Lemma}[section]
\newtheorem{cor}{Corollary}[section]

\newtheorem{rem}{Remark}[section]

\def\H{{\widetilde{H}}}
\def\L{{\Lambda}}
\def\T{{T_0}}
\def\Z{{\Bbb Z}}
\def\Q{{\Bbb Q}}

\def\k{{\overline{k}}}
\def\m{{\mathfrak m}}
\def\p{{\phi}}
\def\s{{\sigma}}

\begin{document}
\title{Maximal Tori determining the Algebraic Group}
\author{Shripad M. Garge}
\address{Harish-Chandra Research Institute, Chhatnag Road, Jhusi, Allahabad. India.}
\email{shripad@mri.ernet.in}

\begin{abstract} 
Let $k$ be a finite field, a global field or a local non-archimedean field. 
Let $H_1$ and $H_2$ be two split, connected, semisimple algebraic groups defined over $k$. 
We prove that if $H_1$ and $H_2$ share the same set of maximal $k$-tori up to $k$-isomorphism, then the Weyl groups $W(H_1)$ and $W(H_2)$ are isomorphic, and hence the algebraic groups modulo their centers are isomorphic except for a switch of a certain number of factors of type $B_n$ and $C_n$. 

We remark that due to a recent result of Philippe Gille, above result holds for fields which admit arbitrary cyclic extensions. 
\end{abstract}  

\maketitle

\section{Introduction.} 
Let $H$ be a connected, semisimple algebraic group defined over a field $k$. 
It is natural to ask to what extent the group $H$ is determined by the $k$-isomorphism classes of maximal $k$-tori contained in it. 
We study this question over finite fields, global fields and local non-archimedean fields.
In this paper, we prove the following theorem. 

\begin{thm}\label{thm:main} $(${\bf Theorem 4.1}$)$
Let $k$ be either a finite field, a global field or a local non-archimedean field. 
Let $H_1$ and $H_2$ be two split, connected, semisimple algebraic groups defined over $k$. 
Suppose that for every maximal $k$-torus $T_1 \subset H_1$ there exists a maximal $k$-torus $T_2 \subset H_2$, such that the tori $T_1$ and $T_2$ are $k$-isomorphic and vice versa. 
Then the Weyl groups $W(H_1)$ and $W(H_2)$ are isomorphic. 

Moreover, if we write the Weyl groups $W(H_1)$ and $W(H_2)$ as a direct product of the Weyl groups of simple algebraic groups, $W(H_1) = \prod_{\L_1} W_{1, \alpha}$, and $W(H_2) = \prod_{\L_2} W_{2, \beta}$, then there exists a bijection $i: \L_1 \rightarrow \L_2$ such that $W_{1, \alpha}$ is isomorphic to $W_{2, i(\alpha)}$ for every $\alpha \in \L_1$.
\end{thm}

Since a split simple algebraic group with trivial center is determined by its Weyl group, except for the groups of the type $B_n$ and $C_n$, we have following theorem. 

\begin{thm}\label{thm:simple}
Let $k$ be as in the previous theorem. 
Let $H_1$ and $H_2$ be two split, connected, semisimple algebraic groups defined over $k$ with trivial center. 
Write $H_i$ as a direct product of simple groups, $H_1 = \prod_{\L_1} H_{1, \alpha}$, and $H_2 = \prod_{\L_2} H_{2, \beta}$. 
If the groups $H_1$ and $H_2$ satisfy the condition given in the above theorem, then there is a bijection $i: \L_1 \rightarrow \L_2$ such that $H_{1, \alpha}$ is isomorphic to $H_{2, i(\alpha)}$, except for the case when $H_{1, \alpha}$ is a simple group of type $B_n$ or $C_n$, in which case $H_{2, i(\alpha)}$ could be of type $C_n$ or $B_n$.
\end{thm}

From the explicit description of maximal $k$-tori in $SO(2n+1)$ and $Sp(2n)$, see for instance \cite[Proposition 2]{K}, one finds that the groups $SO(2n+1)$ and $Sp(2n)$ contain the same set of $k$-isomorphism classes of maximal $k$-tori. 

We show by an example that the existence of split tori in the groups $H_i$ is necessary. 
Note that if $k$ is $\Q_p$, the Brauer group of $k$ is $\Q/\Z$. 
Consider the central division algebras of degree five, $D_1$ and $D_2$, corresponding to $1/5$ and $2/5$ in $\Q/\Z$ respectively, and let $H_1 = SL_1(D_1)$ and $H_2 = SL_1(D_2)$. 
The maximal tori in $H_i$ correspond to the maximal commutative subfields in $D_i$. 
But over $\Q_p$, every division algebra of a fixed degree contains every field extension of that fixed degree (\cite[Proposition 17.10 and Corollary 13.3]{P}), so $H_1$ and $H_2$ share the same set of maximal tori over $k$. 
But they are not isomorphic, since it is known that $SL_1(D) \cong SL_1(D')$ if and only if $D \cong D'$ or $D \cong (D')^{op}$ (\cite[26.11]{KMRT}). 

This paper is arranged as follows. 
The description of the $k$-conjugacy classes of maximal $k$-tori in an algebraic group $H$ defined over $k$ can be given in terms of the Galois cohomology of the normalizer in $H$ of a fixed maximal torus. 
Similarly the $k$-isomorphism classes of $n$-dimensional tori defined over $k$ can be described in terms of $n$-dimensional integral representations of the Galois group of $\k$ (the algebraic closure of $k$) over $k$. 
Using these two descriptions, we obtain a Galois cohomological description for the $k$-isomorphism classes of maximal $k$-tori in $H$, in section 2. 
Since we are dealing with groups that are split over the base field $k$, the Galois action on the Weyl groups is trivial. 
This enables us to prove in section 4, that if the split, connected, semisimple algebraic groups of rank $n$, $H_1$ and $H_2$ share the same set of maximal $k$-tori up to $k$-isomorphism, then the Weyl groups $W(H_1)$ and $W(H_2)$, considered as subgroups of $GL_n(\Z)$, share the same set of elements up to conjugacy in $GL_n(\Z)$. 

This then is the main question to be answered: if Weyl groups of two split, connected, semisimple algebraic groups, $W_1$ and $W_2$, embedded in $GL_n(\Z)$ in the natural way, i.e., by their action on the character group of a fixed split maximal torus, have the property that every element of $W_1$ is $GL_n(\Z)$-conjugate to one in $W_2$, and vice versa, then the Weyl groups are isomorphic. 
Much of the work in this paper is to prove this statement by using elaborate information available about the conjugacy classes in Weyl groups of simple algebraic groups together with their standard representations in $GL_n(\Z)$. 
Our analysis finally depends on the knowledge of characteristic polynomials of elements in the Weyl groups considered as subgroups of $GL_n(\Z)$. 
This information is summarized in section 3. 
Using this information we prove the main theorem in section 4. 

We would like to emphasize that if we were proving Theorems \ref{thm:main} and \ref{thm:simple} for simple algebraic groups, then our proofs are relatively very simple. 
However, for semisimple groups, we have to make a somewhat complicated inductive argument on the maximal rank among the simple factors of the semisimple groups $H_i$. 

We use the term ``simple Weyl group of rank $r$'' for the Weyl group of a simple algebraic group of rank $r$. 
Any Weyl group is a product of simple Weyl groups in a unique way (up to permutation). 
We say that two Weyl groups are isomorphic if and only if the simple factors and their multiplicities are the same. 

The question studied in this paper seems relevant for the study of Mumford-Tate groups over number fields. 
The author was informed, after the completion of the paper, that the Theorem \ref{thm:main} over a finite field is implicit in the work of Larsen and Pink (\cite{LP}). 
We would like to mention that although much of the paper could be said to be implicitly contained in \cite{LP}, the theorems we state (and prove) are not explicitly stated there; besides, our proofs are also different. 

The author wishes to register his acknowledgments towards Prof. Dipendra Prasad for suggesting this question. 
The author enjoyed numerous fruitful discussions with him. 
He spent a lot of time in going through the paper several times and corrected many mistakes. 
The author would like to thank Prof. Gopal Prasad for pointing out a mistake and Prof. J-P. Serre for several useful suggestions. 
He also thanks Prof. M. S. Raghunathan, Prof. R. Parimala and Dr. Maneesh Thakur for encouraging comments and suggestions. 

\section{Galois Cohomological Lemmas.} 
We begin by fixing the notations. 
Let $k$ denote an arbitrary field and $G(\k/k)$ the Galois group of $\k$ (the algebraic closure of $k$) over $k$. 
Let $H$ denote a split, connected, semisimple algebraic group defined over $k$ and let $\T$ be a fixed split maximal torus in $H$. 
Suppose that the dimension of $\T$ is $n$. 
Let $W$ be the Weyl group of $H$ with respect to $\T$. 
Then we have an exact sequence of algebraic groups defined over $k$, 
$$0 ~\longrightarrow ~\T ~\longrightarrow ~N(\T) ~\longrightarrow ~W ~\longrightarrow ~1$$ 
where $N(\T)$ denotes the normalizer of $\T$ in $H$. 

The above exact sequence gives us a map $\psi: H^1(k, N(\T)) \rightarrow H^1(k, W)$. 
It is well known that a certain subset of $H^1(k, N(\T))$ classifies $k$-conjugacy classes of maximal $k$-tori in $H$. 
For the sake of completeness, we formulate it as a lemma in the case of split, connected, semisimple groups, although it is true for a more general class of fields.

\begin{lem}\label{lem:tori:conj}
Let $H$ be a split, connected, semisimple algebraic group defined over a field $k$ and let $\T$ be a fixed split maximal torus in $H$. 
The natural embedding $N(\T) \hookrightarrow H$ induces a map $\Psi: H^1(k, N(\T)) \rightarrow H^1(k, H)$. 
The set of $k$-conjugacy classes of maximal tori in $H$ are in one-one correspondence with the subset of $H^1(k, N(\T))$ which is mapped to the neutral element in $H^1(k, H)$ under the map $\Psi$.
\end{lem}

\begin{proof}
Let $T$ be a maximal $k$-torus in $H$ and let $L$ be a splitting field of $T$, i.e., assume that the torus $T$ splits as a product of ${\Bbb G}_m$'s over $L$. 
We assume that the field $L$ is Galois over $k$. 
By the uniqueness of maximal split tori up to conjugacy, there exists an element $a \in H(L)$ such that $a \T a^{-1} = T$, where $\T$ is the split maximal torus in $H$ fixed before. 
Then for any $\s \in G(L/k)$, we have $\s(a) \T \s(a)^{-1} = T$, as both $\T$ and $T$ are defined over $k$. 
This implies that 
$$\big(a^{-1} \s(a)\big) ~\T ~\big(a^{-1} \s(a)\big)^{-1} = ~\T .$$
Therefore $a^{-1}\s(a) \in N(\T)$. 
This enables us to define a map $:G(L/k) \rightarrow N(\T)$ which sends $\s$ to $a^{-1}\s(a)$. 
By composing this map with the natural map $:G(\k/k) \rightarrow G(L/k)$, we get a map $\p_a : G(\k/k) \rightarrow N(\T)$. 
One checks that 
$$\p_a(\s \tau) = \p_a(\s) \s\big(\p_a(\tau)\big) \hskip3mm \forall ~\s, \tau \in G(\k/k) ,$$ 
i.e., the map $\p_a$ is a 1-cocycle. 
If $b \in H(L)$ is another element such that $b \T b^{-1} = T$, we see that 
$$\p_a(\s) = (b^{-1}a)^{-1} \p_b(\s) \s(b^{-1}a) .$$ 
Thus the element $[\p_a] \in H^1(k, N(\T))$ is determined by the maximal torus $T$. 
We denote it by $\p(T)$. 
It is clear that $\p(T)$ is determined by the $k$-conjugacy class of $T$. 
Moreover, if $\p(T) = \p(S)$ for two maximal tori $T$ and $S$ in $H$, then one can check that these two tori are conjugate over $k$. 
Indeed, if $T = a \T a^{-1}$ and $S = b \T b^{-1}$ for $a, b \in H(\k)$, then for any $\s \in G(\k/k)$, 
$$a^{-1} \s(a) = c^{-1} ~\big(b^{-1} \s(b)\big) ~\s(c)$$ 
for some $c \in N(\T)$. 
Then, $\s(bca^{-1}) = bca^{-1}$ for all $\s \in G(\k/k)$, hence $bca^{-1} \in H(k)$ and $(bca^{-1}) T (bca^{-1})^{-1} = S$. 
Further, it is clear that the image of $\p$ in $H^1(k, N(\T))$ is mapped to the neutral element in $H^1(k, H)$ under the map $\Psi$. 

Moreover, if a 1-cocycle $\p_1 : G(\k/k) \rightarrow N(\T)$ is such that $\Psi(\p_1)$ is neutral in $H^1(k, H)$, then $\phi_1(\s) = a^{-1} \s(a)$ for some $a \in H(\k)$. 
Then the cohomology class $[\p_1] \in H^1(k, N(\T))$ corresponds to the maximal torus $S_1 = a \T a^{-1}$ in $H$. 
Since $a^{-1} \s(a) = \p_1(\s) \in N(\T)$, the torus $S_1$ is invariant under the Galois action, therefore we conclude that it is defined over $k$. 
Thus the image of $\p$ is the inverse image of the neutral element in $H^1(k, H)$ under the map $\Psi$. 
This is the complete description of the $k$-conjugacy classes of maximal $k$-tori in the group $H$. 

Finally, we observe that the detailed proof we have given above amounts to looking at the exact sequence $1 \rightarrow N(T_0) \rightarrow H \rightarrow H/N(T_0) \rightarrow 1$ which gives an exact sequence of pointed sets:
$$H/N(T_0) (k) \longrightarrow H^1(k, N(T_0)) \longrightarrow H^1(k, H) .$$
Therefore $H/N(T_0)(k)$, which is the variety of conjugacy classes of $k$-tori in $H$, is identified to elements in $H^1(k, N(T_0))$ which become trivial in $H^1(k, H)$. 
\end{proof}

We also recall the correspondence between $k$-isomorphism classes of $n$-dimensional $k$-tori and equivalence classes of $n$-dimensional integral representations of $G(\k/k)$. 
Let $\T = {\Bbb G}_m^n$ be the split torus of dimension $n$. 
Let $T_1$ be an $n$-dimensional torus defined over $k$ and let $L_1$ denote the splitting field of $T_1$. 
Since the torus $T_1$ is split over $L_1$, we have an $L_1$-isomorphism $f: \T \rightarrow T_1$. 
The Galois action on $\T$ and $T_1$ gives us another isomorphism, $f^{\s} := \s f \s^{-1} : \T \rightarrow T_1$. 
Again one sees that the map $\varphi_f: G(\k/k) \rightarrow {\rm Aut}_{L_1}(\T)$ given by $\s \mapsto f^{-1}f^{\s}$, is a 1-cocycle. 
Since the torus $\T$ is already split over $k$, we have ${\rm Aut}_{L_1}(\T) \cong {\rm Aut}_k(\T)$, and hence the Galois group $G(\k/k)$ acts trivially on ${\rm Aut}_{L_1}(\T)$, which is isomorphic to $GL_n(\Bbb Z)$. 
Therefore, $\varphi_f$ is actually a homomorphism from the Galois group $G(\k/k)$ to $GL_n(\Bbb Z)$. 
This homomorphism gives an $n$-dimensional integral representation of the absolute Galois group, $G(\k/k)$. 
By changing the isomorphism $f$ to any other $L_1$-isomorphism from $\T$ to $T_1$, we get a conjugate of $\varphi_f$. 
Thus the element $[\varphi_f]$ in $H^1(k, GL_n(\Z))$ is determined by $T_1$ and we denote it by $\varphi(T_1)$. 
Thus a $k$-isomorphism class of an $n$-dimensional torus gives us an equivalence class of $n$-dimensional integral representations of the Galois group, $G(\k/k)$. 
This correspondence is known to be bijective (\cite[2.2]{PR}). 

Since the group $H$ that we consider here is split over the base field $k$, the Weyl group of $H$, which we denote by $W$, is defined over $k$, and $W(\k) = W(k)$. 
Therefore $G(\k/k)$ acts trivially on $W$, and hence $H^1(k, W)$ is the set of conjugacy classes of elements in ${\rm Hom}(G(\k/k), W)$. 
Since $W$ acts faithfully on the character group of $\T$, we can consider $W \hookrightarrow GL_n(\Bbb Z)$ and thus each element of $H^1(k, W)$ gives us an integral representation of the absolute Galois group. 
For a maximal torus $T$ in $H$, we already have an $n$-dimensional integral representation of $G(\k/k)$, as described above. 
We prove that this representation is equivalent to a Galois representation given by an element of $H^1(k, W)$. 

\begin{lem}\label{lem:tori:iso}
Let $H$ be a split, connected, semisimple algebraic group defined over $k$. 
Fix a maximal split $k$-torus $\T$ in $H$. 
Let $T$ be a maximal $k$-torus in $H$. 
Let $\p(T) \in H^1(k, N(\T))$ be the cohomology class corresponding to the $k$-conjugacy class of the torus $T$ in $H$ and $\varphi(T) \in H^1(k, GL_n(\Z))$ be the cohomology class corresponding to the $k$-isomorphism class of $T$. 
Then the integral representations given by $i \circ \psi \circ \p(T)$ and $\varphi(T)$ are equivalent, where $\psi: H^1(k, N(\T)) \rightarrow H^1(k, W)$ is induced by the natural map from $N(\T)$ to $W$, and $i$ is the natural map from $H^1(k, W)$ to $H^1(k, GL_n(\Z))$.
\end{lem}

\begin{proof}
Let $L$ be a splitting field of $T$, then an element $a \in H(L)$ such that $a \T a^{-1} = T$ enables us to define a 1-cocycle $\p_a: G(\k/k) \rightarrow N(\T)$ given by $\p_a(\s) = a^{-1} \s(a)$. 
The element $\p(T) \in H^1(k, N(\T))$ is precisely the element $[\p_a]$.

Further, we treat the conjugation by $a$ as an $L$-isomorphism $f: \T \rightarrow T$, and then it can be checked that the map $f^{\s} := \s f \s^{-1}$ is nothing but the conjugation by $\s(a)$. 
The element $\varphi(T) \in H^1(k, GL_n(\Z))$ is equal to $[\varphi_f]$, where $\varphi_f(\s) = f^{-1}f^{\s}$.
Now, the map $\psi: N(\T) \rightarrow W$ is the natural map taking an element $\alpha \in N(\T)$ to $\overline{\alpha} := \alpha \cdot \T \in W = N(\T)/\T$.
Hence we have 
$$\psi\big(\p_a(\s)\big) = \overline{a^{-1}\s(a)} = f^{-1}f^{\s} = \varphi_f(\s) .$$
Since the action of $W$ on $\T$ is given by conjugation, it is clear that the integral representation of the Galois group, $G(\k/k)$, given by $\psi(\p(T))$ is equivalent to the one given by $\varphi(T)$. 
This proves the lemma.
\end{proof}

Thus, a $k$-isomorphism class of a maximal torus in $H$ gives an element in $H^1(k, W)$. 
We note here that every subgroup of the Weyl group $W$ may not appear as a Galois group of some finite extension $K/k$. 
For instance, if $k$ is a local field of characteristic zero, then it is known that the Galois group of any finite extension over $k$ is a solvable group (\cite[IV]{S}), and thus it is clear that every subgroup of $W$ may not appear as a Galois group. 

If we assume that the base field $k$ is either a finite field or a local non-archimedean field, we have the following result. 

\begin{lem}\label{lem:cyclic}
Let $k$ be a finite field or a local non-archimedean field and let $H$ be a split, connected, semisimple algebraic group defined over $k$. 
Fix a split maximal torus $\T$ in $H$ and let $W$ denote the Weyl group of $H$ with respect to $\T$. 
An element in $H^1(k, W)$ which corresponds to a homomorphism $\rho: G(\k/k) \rightarrow W$ with cyclic image, corresponds to a $k$-isomorphism class of a maximal torus in $H$ under the mapping $\psi: H^1(k, N(\T)) \rightarrow H^1(k, W)$.
\end{lem}

\begin{proof}
Consider the map $\Psi:H^1(k, N(\T)) \rightarrow H^1(k, H)$ induced by the inclusion $N(\T) \hookrightarrow H$. 
If we denote the neutral element in $H^1(k, H)$ by $\iota$, then by Lemma \ref{lem:tori:conj} the set 
$$X := \big\{f \in H^1(k, N(\T)): \Psi(f) = \iota\big\}$$ 
is in one-one correspondence with the $k$-conjugacy classes of maximal $k$-tori in $H$. 
By the previous lemma, it is enough to show that $[\rho] \in \psi(X)$, where $\psi: H^1(k, N(\T)) \rightarrow H^1(k, W)$ is induced by the natural map from $N(\T)$ to $W$. 

By Tits' theorem (\cite[4.6]{T}), there exists a subgroup $\overline{W}$ of $N(\T)(k)$ such that the sequence 
$$0 ~\longrightarrow ~\mu_2^n ~\longrightarrow ~\overline{W} ~\longrightarrow ~W ~\longrightarrow ~1$$
is exact. 
Let $N$ denote the image of $\rho$ in $W$. 
We know that $N$ is a cyclic subgroup of $W$.
Let $w$ be a generator of $N$ and $\overline{w}$ be a lifting of $w$ to $\overline{W}$. 
Since the base field $k$ admits cyclic extensions of any given degree, there exists a map $\rho_1$ from $G(\k/k)$ to $\overline{W}$ whose image is the cyclic subgroup generated by $\overline{w}$. 
Since the Galois action on $\overline{W}$ is trivial, as $\overline{W}$ is a subgroup of $N(\T)(k)$, the map $\rho_1$ could be treated as a 1-cocycle from $G(\k/k)$ to $N(\T)$. 
Consider $[\rho_1]$ as an element in $H^1(k, N(\T))$, then $\psi([\rho_1]) = [\rho] \in H^1(k, W)$. 
Now we make two cases.

\vskip2mm
\noindent{Case 1: $k$ is a finite field.} 

By Lang's Theorem (\cite[Corollary to Theorem 1]{L}), $H^1(k, H)$ is trivial and hence the set $X$ coincides with $H^1(k, N(\T))$. 
Therefore the element $[\rho_1] \in H^1(k, N(\T))$ corresponds to a $k$-conjugacy class of maximal $k$-torus in $H$.
Then, by previous lemma, $[\rho] = \psi([\rho_1])$ corresponds to a $k$-isomorphism class of maximal $k$-tori in $H$. 

\vskip2mm
\noindent{Case 2: $k$ is a local non-archimedean field.} 

By \cite[Proposition 2.10]{PR}, there exists a semisimple simply connected algebraic group $\H$, which is defined over $k$, together with a $k$-isogeny $\pi: \H \rightarrow H$. 
We have already fixed a split maximal torus $\T$ in $H$, let $\widetilde{\T}$ be the split maximal torus in $\H$ which gets mapped to $\T$ by the covering map $\pi$. 
It can be seen that, by restriction we get a surjective map $\pi: N(\widetilde{\T}) \rightarrow N(\T)$, where the normalizers are taken in appropriate groups. 
Moreover, the induced map $\pi_1: \widetilde{W} \rightarrow W$ is an isomorphism.

We define the maps $$\widetilde{\psi}: H^1(k, N(\widetilde{\T})) \rightarrow H^1(k, \widetilde{W})$$ and $$\widetilde{\Psi}: H^1(k, N(\widetilde{\T})) \rightarrow H^1(k, \H)$$ in the same way as the maps $\psi$ and $\Psi$ are defined for the group $H$.

Consider the following diagram, where the horizontal arrows represent natural maps. 
$$
\begin{CD}
  \widetilde{H} @<<< N(\widetilde{T_0}) @>>> \widetilde{W} \\
  @V{\pi}VV   @V{\pi}VV @VV{\pi_1}V \\
  H @<<< N(T_0) @>>> W ,
\end{CD}
$$
It is clear that the above diagram is commutative and hence so is the following diagram. 
$$
\begin{CD}
  H^1(k, \widetilde{H}) @<{\widetilde{\Psi}}<< H^1(k, N(\widetilde{T_0})) @>{\widetilde{\psi}}>> H^1(k, \widetilde{W}) \\
  @V{\pi^*}VV   @V{\pi^*}VV @VV{\pi_1^*}V \\
  H^1(k, H) @<<{\Psi}< H^1(k, N(T_0)) @>>{\psi}> H^1(k, W) .
\end{CD}
$$
Since $\pi_1$ is an isomorphism, the map $\pi_1^*$ is a bijection. 
Now, consider an element $[\rho] \in H^1(k, W)$, such that the image of the $1$-cocycle $\rho$ is a cyclic subgroup of $W$, and let $[\widetilde{\rho}]$ be its inverse image in $H^1(k, \widetilde{W})$ under the bijection $\pi_1^*$. 
Using Tits' theorem (\cite{T}) as above, we lift $[\widetilde{\rho}]$ to an element $[\widetilde{\rho_1}]$ in $H^1(k, N(\widetilde{\T}))$. 
Since $\H$ is simply connected and $k$ is a non-archimedean local field, $H^1(k, \H)$ is trivial (\cite{BT}, \cite{Kn1}, \cite{Kn2}). 
Therefore, $\widetilde{\Psi}([\widetilde{\rho_1}])$ is neutral in $H^1(k, \H)$ and so is $\pi^*(\widetilde{\Psi}([\widetilde{\rho_1}]))$ in $H^1(k, H)$.
By commutativity of the diagram, we have that the element $[\rho] \in H^1(k, W)$ has a lift $\pi^*([\widetilde{\rho_1}])$ in $H^1(k, N(\T))$ such that $\Psi(\pi^*([\widetilde{\rho_1}]))$ is neutral in $H^1(k, H)$. 
Thus the element $[\rho]$ corresponds to a $k$-isomorphism class of a maximal torus in $H$. 

This proves the lemma.
\end{proof}

\section{Characteristic Polynomials.}

For a finite subgroup $W$ of $GL_n(\Z)$, we define $ch(W)$ to be the set of characteristic polynomials of elements of $W$ and $ch^*(W)$ to be the set of irreducible factors of elements of $ch(W)$.
Since all the elements of $W$ are of finite order, the irreducible factors (over $\Q$) of the characteristic polynomials are cyclotomic polynomials. 
We denote by $\p_r$, the $r$-th cyclotomic polynomial, i.e., the irreducible monic polynomial over $\Z$ satisfied by a primitive $r$-th root of unity. 
We define 
$$\m_i(W) = \max \big\{t: \p_i^t {\rm ~divides} ~f {\rm ~for ~some} ~f \in ch(W)\big\}$$
and
$$\m_i'(W) = \min \big\{t: \p_2^t \cdot \p_i^{\m_i(W)} {\rm ~divides} ~f {\rm ~for ~some} ~f \in ch(W)\big\} .$$
For positive integers $i \ne j$, we define 
$$\m_{i, j}(W) = \max \big\{t+s: \p_i^t \cdot \p_j^s {\rm ~divides} ~f {\rm ~for ~some} ~f \in ch(W)\big\} .$$

If $U_1$ is a subgroup of $GL_n(\Z)$ and $U_2$ is a subgroup of $GL_m(\Z)$, then $U_1 \times U_2$ can be treated as a subgroup of $GL_{m+n}(\Z)$.
Then
$$ch(U_1 \times U_2) = \big\{f_1 \cdot f_2: f_1 \in ch(U_1), f_2 \in ch(U_2)\big\} .$$
Moreover, one can easily check that 
$$\m_i(U_1 \times U_2) = \m_i(U_1) + \m_i(U_2) \hskip3mm \forall ~i ,$$ 
$$\m_i'(U_1 \times U_2) = \m_i'(U_1) + \m_i'(U_2) \hskip3mm \forall ~i$$ 
and 
$$\m_{i, j}(U_1 \times U_2) = \m_{i, j}(U_1) + \m_{i, j}(U_2) \hskip3mm \forall ~i, j .$$ 
A simple Weyl group $W$ of rank $n$ has a natural embedding in $GL_n(\Z)$. 
We obtain a description of the sets $ch^*(W)$ with respect to this natural embedding. 
Here we use the following result due to T. A. Springer (\cite[Theorem 3.4(i)]{Sp}) about the fundamental degrees of the Weyl group $W$. 
We recall that the degrees of the generators of the invariant algebra of the Weyl group are called as the fundamental degrees of the Weyl group. 

\begin{thm}[Springer]\label{thm:Sp}
Let $W$ be a complex reflection group and $d_1, d_2, \dots, d_m$ be the fundamental degrees of $W$. 
An $r^{th}$ root of unity occurs as an eigenvalue for some element of $W$ if and only if $r$ divides one of the fundamental degrees $d_i$ of $W$. 

Equivalently, the irreducible polynomial $\p_r$ is in $ch^*(W)$ if and only if $r$ divides one of the fundamental degrees $d_i$ of the reflection group $W$.
\end{thm}

We now list the fundamental degrees and the divisors of degrees for the simple Weyl groups, cf. \cite[3.7]{H}:

\newpage
\begin{center}
{\bf Table 3.2.}
\vskip2mm
\begin{tabular}{|c|l|l|}
\hline
Type & Degrees & Divisors of degrees \\ \hline
$A_n$ & $2, 3, \dots, n + 1$ & $1, 2, \dots, n + 1$ \\
$B_n$ & $2, 4, \dots, 2n$ & $1, 2, \dots, n, n + 2, n + 4, \dots, 2n$ \hfill $n$ even\\
 & & $1, 2, \dots, n, n + 1, n + 3, \dots, 2n$ \hfill $n$ odd\\
$D_n$ & $2, 4, \dots, 2n - 2, n$ & $1, 2, \dots, n, n + 2, n + 4, \dots, 2n - 2$ \hfill $n$ even\\
 & & $1, 2, \dots, n, n + 1, n + 3, \dots, 2n - 2$ \hfill $n$ odd\\
$G_2$ & $2, 6$ & $1, 2, 3, 6$ \\
$F_4$ & $2, 6, 8, 12$ & $1, 2, 3, 4, 6, 8, 12$ \\
$E_6$ & $2, 5, 6, 8, 9, 12$ & $1, 2, 3, 4, 5, 6, 8, 9, 12$ \\
$E_7$ & $2, 6, 8, 10, 12, 14, 18$ & $1, 2, 3, 4, 5, 6, 7, 8, 9, 10, 12, 14, 18$ \\
$E_8$ & $2, 8, 12, 14, 18, 20, 24, 30$ & $1, 2, 3, 4, 5, 6, 7, 8, 9, 10, 12, 14, 15, 18, 20, 24, 30$ \\
\hline
\end{tabular}
\end{center}
\vskip2mm
Using Theorem \ref{thm:Sp} and the above table, we can now easily compute the set $ch^*(W)$ for any simple Weyl group $W$. 
We summarize them below. 
\begin{eqnarray}
ch^*\big(W(A_n)\big) & = & \big\{\p_1, \p_2, \dots, \p_{n + 1}\big\} \\
ch^*\big(W(B_n)\big) & = & \big\{\p_i, \p_{2i}: i = 1, 2, \dots, n\big\} \\
ch^*\big(W(D_n)\big) & = & \big\{\p_i, \p_{2j}: i = 1, 2, \dots, n, j = 1, 2 \dots, n - 1 \big\} \\
ch^*\big(W(G_2)\big) & = & \big\{\p_1, \p_2, \p_3, \p_6\big\} \\
ch^*\big(W(F_4)\big) & = & \big\{\p_1, \p_2, \p_3, \p_4, \p_6, \p_8, \p_{12}\big\} \\
ch^*\big(W(E_6)\big) & = & \big\{\p_1, \p_2, \p_3, \p_4, \p_5, \p_6, \p_8, \p_9, \p_{12}\big\} \\
ch^*\big(W(E_7)\big) & = & \big\{\p_1, \p_2, \dots, \p_{10}, \p_{12}, \p_{14}, \p_{18}\big\} \\
ch^*\big(W(E_8)\big) & = & \big\{\p_1, \p_2, \dots, \p_{10}, \p_{12}, \p_{14}, \p_{15}, \p_{18}, \p_{20}, \p_{24}, \p_{30}\big\}
\end{eqnarray}

\section{Main Result.} 

In this section, $k$ is either a finite field, a global field or a non-archimedean local field. 
We now restate the main theorem, the Theorem \ref{thm:main} of the introduction. 

\begin{thm}\label{thm:main2}
Let $H_1$ and $H_2$ be two split, connected, semisimple algebraic groups defined over $k$. 
Suppose that for every maximal $k$-torus $T_1 \subset H_1$ there exists a maximal $k$-torus $T_2 \subset H_2$ such that the torus $T_2$ is $k$-isomorphic to the torus $T_1$ and vice versa. 
Then, the Weyl groups $W(H_1)$ and $W(H_2)$ are isomorphic. 

Moreover, if we write the Weyl groups $W(H_1)$ and $W(H_2)$ as a direct product of the Weyl groups of simple algebraic groups, $W(H_1) = \prod_{\L_1} W_{1, \alpha}$, and $W(H_2) = \prod_{\L_2} W_{2, \beta}$, then there exists a bijection $i: \L_1 \rightarrow \L_2$ such that $W_{1, \alpha}$ is isomorphic to $W_{2, i(\alpha)}$ for every $\alpha \in \L_1$.
\end{thm}

The proof of this theorem occupies the rest of this section. 
Clearly, the above mentioned groups $H_1$ and $H_2$ are of the same rank, say $n$. 
Let $W_1$ and $W_2$ denote the Weyl groups of $H_1$ and $H_2$, respectively. 
We always treat the Weyl groups $W_1$ and $W_2$ as subgroups of $GL_n(\Z)$. 
We first prove a lemma here which transforms the information about $k$-isomorphism of maximal $k$-tori in the groups $H_i$ into some information about the conjugacy classes of the elements of the corresponding Weyl groups $W_i$. 

\begin{lem}\label{lem:Weyl}
Under the hypothesis of theorem $\ref{thm:main2}$, for every element $w_1 \in W_1$ there exists an element $w_2 \in W_2$ such that $w_2$ is conjugate to $w_1$ in $GL_n(\Z)$ and vice versa.
\end{lem}

\begin{proof}
Let $w_1 \in W_1$ and let $N_1$ denote the subgroup of $W_1$ generated by $w_1$. 
Since the base field $k$ admits any cyclic group as a Galois group, there is a map $\rho_1 : G(\k/k) \rightarrow W_1$ such that $\rho_1(G(\k/k)) = N_1$. 

We first consider the case when $k$ is a finite field or a local non-archimedean field. 
By Lemma \ref{lem:cyclic}, the element $[\rho_1] \in H^1(k, W_1)$ corresponds to a maximal $k$-torus in $H_1$, say $T_1$. 
By the hypothesis, there exists a torus $T_2 \subset H_2$ which is $k$-isomorphic to $T_1$. 
We know by Lemma \ref{lem:tori:iso} that there exists an integral Galois representation $\rho_2 : G(\k/k) \rightarrow GL_n(\Z)$ corresponding to the $k$-isomorphism class of $T_2$ which factors through $W_2$. 
Let $N_2 := \rho_2(G(\k/k)) \subseteq W_2$. 
Since $T_1$ and $T_2$ are $k$-isomorphic tori, the corresponding Galois representations, $\rho_1$ and $\rho_2$, are equivalent. 
This implies that there exists $g \in GL_n(\Z)$ such that $N_2 = g N_1 g^{-1}$. 
Then $w_2 := g w_1 g^{-1} \in N_2 \subseteq W_2$ is a conjugate of $w_1$ in $GL_n(\Z)$. 
We can start with an element $w_2 \in W_2$ and obtain its $GL_n(\Z)$-conjugate in $W_1$ in the same way. 

Now we consider the case when $k$ is a global field. 
Let $v$ be a non-archimedean valuation of $k$ and let $k_v$ be the completion of $k$ with respect to $v$. 
Clearly the groups $H_1$ and $H_2$ are defined over $k_v$. 
Let $T_{1, v}$ be a maximal $k_v$-torus in $H_1$. 
Then by Grothendieck's theorem (\cite[7.9, 7.11]{BS}) and weak approximation property (\cite[Proposition 7.3]{PR}), there exists a $k$-torus in $H$, say $T_1$, such that $T_{1, v}$ is obtained from $T_1$ by the base change. 
By hypothesis, we have a $k$-torus $T_2$ in $H_2$ which is $k$-isomorphic to $T_1$. 
Then the torus $T_{2, v}$, obtained from $T_2$ by the base change, is $k_v$-isomorphic to $T_{1, v}$. 
Thus, every maximal $k_v$-torus in $H_1$ has a $k_v$-isomorphic torus in $H_2$. 
Similarly, we can show that every maximal $k_v$-torus in $H_2$ has a $k_v$-isomorphic torus in $H_1$. 
Then, the proof follows by previous case. 

This proves the lemma.
\end{proof}

\begin{cor}\label{cor}
Under the hypothesis of theorem $\ref{thm:main2}$, $ch(W_1) = ch(W_2)$ and $ch^*(W_1) = ch^*(W_2)$. 
In particular, $\m_i(W_1) = \m_i(W_2)$, $\m_i'(W_1) = \m_i'(W_2)$ and $\m_{i, j}(W_1) = \m_{i, j}(W_2)$ for all $i, j$.
\end{cor}

\begin{proof}
Since the Weyl groups $W_1$ and $W_2$ share the same set of elements up to conjugacy in $GL_n(\Z)$, the sets $ch(W_i)$ are the same for $i = 1, 2$, and hence the sets $ch^*(W_i)$ are also the same for $i = 1, 2$. 
Further, for a fixed integer $i$, $\p_i^{\m_i(W_1)}$ divides an element $f_1 \in ch(W_1)$. 
But since $ch(W_1) = ch(W_2)$, the polynomial $\p_i^{\m_i(W_1)}$ also divides an element $f_2 \in ch(W_2)$. 
Therefore $\m_i(W_1) \leq \m_i(W_2)$. 
We obtain the inequality in the other direction in the same way and hence $\m_i(W_1) = \m_i(W_2)$. 
Similarly, we can prove that $\m_i'(W_1) = \m_i'(W_2)$.

Similarly, for integers $i \ne j$, the sets 
$$\big\{(t_1, s_1): \p_i^{t_1} \cdot \p_j^{s_1} {\rm ~divides ~some ~element} ~f_1 \in ch(W_1)\big\} ,$$ 
$$\big\{(t_2, s_2): \p_i^{t_2} \cdot \p_j^{s_2} {\rm ~divides ~some ~element} ~f_2 \in ch(W_2)\big\}$$ 
are the same for $i = 1, 2$. 
It follows that $\m_{i, j}(W_1) = \m_{i, j}(W_2)$.
\end{proof}

We now prove the following result before going on to prove the main theorem. 

\begin{thm}\label{thm:ind}
Let $H_1$ and $H_2$ be split, connected, semisimple algebraic groups of rank $n$. 
Assume that $\m_i(W(H_1)) = \m_i(W(H_2))$, $\m_i'(W(H_1)) = \m_i'(W(H_2))$ and $\m_{i, j}(W(H_1)) = \m_{i, j}(W(H_2))$ for all $i, j$. 
Let $m$ be the maximum possible rank among the simple factors of $H_1$ and $H_2$. 
Let $W_i'$ denote the product of the Weyl groups of simple factors of $H_i$ of rank $m$ for $i = 1, 2$. 
Then the groups $W_1'$ and $W_2'$ are isomorphic.
\end{thm}

\begin{proof}
We denote $W(H_i)$ by $W_i$ for $i = 1, 2$. 
We prove that if a simple Weyl group of rank $m$ appears as a factor of $W_1$ with multiplicity $p$, then it appears as a factor of $W_2$, with the same multiplicity. 
We prove this lemma case by case depending on the type of rank $m$ simple factors of $H_i$. 

We prove this result by comparing the sets $ch^*(W)$ for the simple Weyl groups of rank $m$. 
We observe from the Table 3.2, that the maximal degree of the simple Weyl group of exceptional type, if any, is the largest among the maximal degrees of simple Weyl groups of rank $m$. 
The next largest maximal degree is that of $W(B_m)$, the next one is that of $W(D_m)$ and finally the Weyl group $W(A_m)$ has the smallest maximal degree. 
We use the relation between the elements of $ch^*(W)$ and the degrees of the Weyl group $W$, given by Theorem \ref{thm:Sp}. 
So, we begin the proof of the lemma with the case of exceptional groups of rank $m$, prove that it occurs with the same multiplicity for $i = 1, 2$. 
Then we prove the lemma for $B_m$, then for $D_m$ and finally we prove the lemma for the group $A_m$. 

\vskip2mm
\noindent{\bf Case 1.} One of the groups $H_i$ contains a simple exceptional factor of rank $m$. 

We first treat the case of the simple group $E_8$, i.e., we assume that 8 is the maximum possible rank of the simple factors of the groups $H_i$. 
We know that $\m_{30}(W(E_8)) = 1$. 
Observe that $\p_{30}$ is an irreducible polynomial of degree 8, hence it cannot occur in $ch^*(W)$ for any simple Weyl group of rank $\leq 7$. 
Moreover, from the Theorem \ref{thm:Sp} and the Table 3.2, it is clear that $\m_{30}(W(A_8)) = \m_{30}(W(B_8)) = \m_{30}(W(D_8)) = 0$. 
Hence the multiplicity of $E_8$ in $H_i$ is given by $\m_{30}(W_i)$ which is the same for $i = 1, 2$. 

Similarly for the simple algebraic group $E_7$, we observe that $\m_{18}(W(E_7)) = 1$ and $\m_{18}(W) = 0$ for any simple Weyl group $W$ of rank $\leq 7$. 
Then the multiplicity of $E_7$ in $H_i$ is given by $\m_{18}(W_i)$ which is the same for $i = 1, 2$. 

The case of $E_6$ is done by using $\m_9$. 
It is clear that $\m_9(W) = 0$ for any simple Weyl group $W$ of rank $\leq 6$. 

The cases of $F_4$ and $G_2$ are done similarly by using $\m_{12}$ and $\m_6$ respectively. 

\vskip2mm
\noindent{\bf Case 2.} One of the groups $H_i$ has $B_m$ or $C_m$ as a factor. 

Since $W(B_m) \cong W(C_m)$, we treat the case of $B_m$ only. 
By case 1, we can assume that the exceptional group of rank $m$, if any, occurs with the same multiplicities in both $H_1$ and $H_2$, and hence while counting the multiplicities $\m_i$, $\m_i'$ and $\m_{i, j}$, we can (and will) ignore the exceptional groups of rank $m$. 

Observe that $\m_{2m}(W(B_m)) = 1$ and $\m_{2m}(W) = 0$ for any other simple Weyl group $W$ of classical type of rank $\leq m$. 
However, it is possible that $\m_{2m}(W) \ne 0$ for a simple Weyl group $W$ of exceptional type of rank strictly less than $m$. 
If $m \geq 16$, then this problem does not arise, therefore the multiplicity of $B_m$ in $H_i$ for $m \geq 16$ is given by $\m_{2m}(W_i)$, which is the same for $i= 1, 2$. 
We do the cases of $B_m$ for $m \leq 15$ separately. 

For the group $B_2$, we observe that $\m_4(W(B_2)) = 1$ and $\m_4(W) = 0$ for any other simple Weyl group $W$ of rank $\leq 2$. 
Thus, the case of $B_2$ is done using $\m_4(W_1) = \m_4(W_2)$. 

For the group $B_3$, we have $\m_6(W(B_3)) = 1$, but then $\m_6(W(G_2))$ is also $1$. 
Observe that $\m_4(W(B_3)) = 1$ and $\m_4(W(G_2)) = 0$. 
We do this case by looking at the multiplicities of $\p_4$ and $\p_6$, so we do not worry about the simple Weyl groups $W$ of rank $\leq 3$ for which the multiplicities $\m_4(W)$ and $\m_6(W)$ are both zero. 
Now, let the multiplicities of $B_3$, $G_2$ and $B_2$ in the groups $H_i$ be $p_i$, $q_i$ and $r_i$, for $i = 1, 2$ respectively. 
Then, using $\m_6(W_1) = \m_6(W_2)$, we get that $p_1 + q_1 = p_2 + q_2$. 
Using $\m_4$ we have that $p_1 + r_1 = p_2 + r_2$ and using $\m_{4, 6}$ we get $p_1 + q_1 + r_1 = p_2 + q_2 + r_2$. 
Combining these equalities, we get that $p_1 = p_2$, i.e., the group $B_3$ appears in both the groups $H_i$ with the same multiplicity. 

For the group $B_4$, we observe that $\m_8(W(B_4)) = 1$. 
Since $\p_8$ has degree $4$, it cannot occur in $ch(W)$ for any simple Weyl group of rank $\leq 3$ and $\m_8(W(A_4)) = \m_8(W(D_4)) = 0$. 
Since we are assuming by case $1$ that the group $F_4$ occurs in both $H_i$ with the same multiplicity, we are done in this case also. 

For the group $B_5$, we have $\m_{10}(W(B_5)) = 1$ and $\m_{10}(W) = 0$ for any other simple Weyl group of classical type of rank $\leq 5$. 
Since $5$ does not divide the order of $W(G_2)$ or $W(F_4)$, $\m_{10}(W(G_2)) = \m_{10}(W(F_4)) = 0$ and so we are done. 

The group $B_6$ is another group where the exceptional groups give problems. 
We have $\m_{12}(W(B_6)) = 1$, but $\m_{12}(W(F_4))$ is also $1$. 
Observe that $\m_{10}(W(B_6)) = 1$, but $\m_{10}(W(F_4)) = 0$. 
Now, let the multiplicities of $B_6, D_6, B_5$ and $F_4$ in $H_i$ be $p_i, q_i, r_i$ and $s_i$ respectively. 
Then, 
$$p_1 + s_1 = \m_{12}(W_1) = \m_{12}(W_2) = p_2 + s_2 .$$ 
Similarly comparing $\m_{10}$, we get that 
$$p_1 + q_1 + r_1 = p_2 + q_2 + r_2 .$$ 
Then, we compare $\m_{10, 12}$ of the groups $W_1$ and $W_2$, to get that 
$$p_1 + q_1 + r_1 + s_1 = p_2 + q_2 + r_2 + s_2 .$$ 
Combining this equality with the one obtained by $\m_{10}$, we get that $s_1 = s_2$ and hence $p_1 = p_2$. 
Thus the group $B_6$ occurs in both $H_1$ and $H_2$ with the same multiplicity. 

We have that $\m_{14}(W(E_6)) = 0$, therefore the group $B_7$ is characterized by $\p_{14}$ and hence it occurs in both $H_i$ with the same multiplicity. 

For the group $B_8$, $\m_{16}(W(B_8)) = 1$. 
Since $\p_{16}$ has degree $8$, it cannot occur in $ch^*(W)$ for any of the Weyl groups of $G_2, F_4, E_6$ or $E_7$. 
Thus, the group $B_8$ is characterized by $\p_{16}$ and hence it occurs in both $H_i$ with the same multiplicity. 

The group $B_9$ has the property that $\m_{18}(W(B_9)) = 1$. 
But we also have that $\m_{18}(W(E_7)) = \m_{18}(W(E_8)) = 1$. 
We conclude that the multiplicity of $E_8$ is the same for both $W_1$ and $W_2$ using $\m_{30}$. 
Then we compare the multiplicities $\m_{18}, \m_{16}$ and $\m_{16, 18}$ to prove that the group $B_9$ occurs in both the groups $H_i$ with the same multiplicity. 

Now, we do the case of $B_{10}$. 
Here $\m_{20}(W(B_{10})) = 1$. 
Observe that $\m_{20}(W) = 0$, for any other simple Weyl group $W$ of rank $\leq 10$, except for $E_8$. 
Then the multiplicity of $B_{10}$ in $H_i$ is $\m_{20}(W_i) - \m_{30}(W_i)$ and hence it is the same for $i = 1, 2$. 

The same method works for $B_{12}$ also, i.e. the multiplicity of $B_{12}$ in $H_i$ is $\m_{24}(W_i) - \m_{30}(W_i)$. 

The multiplicities of $B_{11}, B_{13}$ and $B_{14}$ in $H_i$ are given by $\m_{22}(W_i), \m_{26}(W_i)$ and $\m_{28}(W_i)$ and hence they are the same for $i = 1, 2$. 
Now we are left with the case of $B_{15}$ only.  

For $B_{15}$, we have $\m_{30}(W(B_{15})) = \m_{30}(W(E_8)) = 1$ and it is $0$ for any other simple Weyl group of rank $\leq 15$. 
Observe that $\m_{28}(W(B_{15})) = \m_{28}(W(B_{14})) = 1$ and it is $0$ for any other simple Weyl group of rank $\leq 15$. 
Then by comparing $\m_{30}$, $\m_{28}$ and $\m_{28, 30}$ we get the desired result that $B_{15}$ occurs in both $H_i$ with the same multiplicity. 

\vskip2mm
\noindent{\bf Case 3.} One of the groups $H_i$ has $D_m$ as a factor. 

While doing the case of $D_m$, we assume that the exceptional group of rank $m$, if any, and the group $B_m$ occur in both $H_i$ with the same multiplicities. 

We observe that $2m - 2$ is the largest integer $r$ such that $\p_r \in ch^*(W(D_m))$, but $\m_{2m -2}(W(B_{m-1})) = 1$. 
Hence we always have to compare the group $D_m$ with the group $B_{m - 1}$. 

Let us assume that $m \geq 17$, so that $\p_{2m -2} \not \in ch^*(W)$ for any simple Weyl group of exceptional type of rank $< m$. 

We know that $\m_{2m-2}(W(D_m)) = \m_{2m - 2}(W(B_{m - 1})) = 1$ and for any other simple Weyl group $W$ of classical type of rank $\leq m$, $\m_{2m - 2}(W) = 0$. 
Further, $(X+1)(X^{m-1}+1)$ is the only element in $ch(W(D_m))$ which has $\p_{2m-2}$ as a factor. 
Similarly $X^{m-1}+1$ is the only element in $ch(W(B_{m-1}))$ which has $\p_{2m-2}$ as a factor. 
Observe that $\m_{2m -2}'(W(D_m)) = \m_{2m -2}'(W(B_{m - 1})) + 1$ and $\ m_{2m -2}'(W) = 0$ for any other simple Weyl group $W$ of rank $\leq m$. 
Now, let $p_i$ and $q_i$ be the multiplicities of the groups $D_m$ and $B_{m-1}$ in $H_i$ for $i=1, 2$ respectively. 
Then by considering $\m_{2m -2}$, we have 
$$p_1 + q_1 = p_2 + q_2 .$$ 
Further if $m$ is even, then by considering $\m_{2m -2}'$ we have 
$$2p_1 + q_1 = 2p_2 + q_2 .$$ 
This equality, combined with the previous equality, implies that $p_1 = p_2$. 
If $m$ is odd then $\m_{2m -2}'$ itself gives that $p_1 = p_2$. 
Thus, we get the result that the group $D_m$ appears in both $H_i$ with the same multiplicity for $i = 1, 2$.

Now we do the cases of $D_m$, for $m \leq 16$. 

For the group $D_4$, we have to consider the simple algebraic groups $B_3$ and $G_2$. 
Comparing the multiplicities $\m_6$, $\m_4$ and $\m_{4, 6}$ we get that $G_2$ occurs in both $H_i$ with the same multiplicity, and then we proceed as above to prove that $D_4$ also occurs with the same multiplicity in both the groups $H_i$. 

For the group $D_5$, we first prove that the multiplicity of $F_4$ is the same for both $H_i$ using $\m_{12}$ and then prove the required result by considering $\m_5$, $\m_8$ and $\m_{5, 8}$. 
Now, while dealing the case of $D_6$, we observe that $\m_{10}(W(G_2)) = \m_{10}(W(F_4)) = 0$, and so we do this case as done above for $m \geq 17$. 
The case of $D_7$ is done by considering $\m_7$, $\m_{12}$ and $\m_{7, 12}$. 
While doing the case of $D_8$, we first prove that the group $E_7$ occurs in both the $H_i$ with the same multiplicity by considering $\m_{18}$ and then proceed as above. 
For the group $D_9$, we prove that $E_8$ occurs in both $H_i$ with the same multiplicity by considering $\m_{30}$ and proceed as done above for $m \geq 17$. 
While doing the case $D_{10}$, we prove that $E_8$ appears in both the $H_i$ with the same multiplicity by considering $\m_{30}$ and the same can be proved for $E_7$ by considering $\m_{18}$, $\m_{16}$ and $\m_{16, 18}$. Then we do this case as done above
. 

For the groups $D_m$, $m \geq 11$, the only simple Weyl group $W$ of exceptional type such that $\p_{2m -2} \in ch^*(W)$ is $W(E_8)$, but for $D_m$, $m \leq 15$, we can assume that $E_8$ occurs in both the $H_i$ with the same multiplicity by considering $\m_{30}$ and hence we are through. 
For the group $D_{16}$, we take care of $E_8$ by considering $\m_{30}$, $\m_{28}$ and $\m_{28, 30}$. 
Other arguments are the same as done in the case $m \geq 17$. 

\vskip2mm
\noindent{\bf Case 4.} One of the groups $H_i$ has $A_m$ as a factor. 

Now we do the last case, the case of simple algebraic group of type $A_m$. 
Here, as usual, we assume that all other simple algebraic groups of rank $m$ occur with the same multiplicities in both $H_i$. 

If $m$ is even, then $m+1$ is odd and hence $\m_{m+1}(W) = 0$, for any simple Weyl group $W$ of classical type of rank $< m$. 
If $m \geq 30$, then we do not have to bother about the exceptional simple groups of rank $< m$. 
If $m$ is odd and $m \geq 31$, then $\p_{m+1}$ occurs in $ch^*(W(B_r))$ and $ch^*(W(D_{r+1}))$ for $r \geq (m+1)/2$. 
Then we compare the multiplicities $\m_m$, $\m_{m+1}$ and $\m_{m, m+1}$ and get the result that the group $A_m$ occurs in $H_i$ with the same multiplicity. 
So, we have to do the case of $A_m$ for $m \leq 29$, separately. 

Now, we consider the cases $A_m$ for $m \leq 29$. 
The cases of $A_1, A_2$ are easy since there are no exceptional groups of rank $1$. 
For $A_3$, we use $\m_3, \m_4$ and $\m_{3, 4}$ to get the result. 
Similarly $A_4$ is done by using $\m_5$. 
The problem comes for $A_5$, since $\m_6(W(B_3)) \ne 0$, $\m_6(W(G_2)) \ne 0$ and $\m_6(W(F_4)) \ne 0$. 
But, this is handled by first proving that $F_4$ appears with the same multiplicity using $\m_{12}$ and then using the multiplicities $\m_5$, $\m_6$ and $\m_{5, 6}$. 
The case of $A_6$ is done by using $\m_7$. 
For $A_7$, we use $\m_7$, $\m_8$ and $\m_{7, 8}$. 
 
While doing the case of $A_8$, we can first assume that the multiplicity of $E_7$ is the same for both $H_i$, by using $\m_{18}$. 
Then we use $\m_7$, $\m_9$ and $\m_{7, 9}$ to get the result. 
For the group $A_9$, we can again get rid of $E_7$ and $E_8$ using the multiplicities $\m_{18}$ and $\m_{30}$. 
Then we are left with the groups $B_5$ and $E_6$, so here we work with $\m_7$, $\m_{10}$ and $\m_{7, 10}$ to get the result. 

Further, we observe that for $m \in \{10, 12, \dots, 28\}$ such that $m \ne 14$, we have that $\m_{m + 1} (W) = 0$ for any simple Weyl group of rank $< m$. 
Thus, the multiplicities of the groups $A_m$, where $m \in \{10, 12, \dots, 28\}$ and $m \ne 14$, in $H_i$ are characterized by considering $\m_{m + 1}(W_i)$ and hence they are same for $i = 1, 2$. 
The case of $A_{14}$ is done by using $\m_{13}, \m_{15}$ and $\m_{13, 15}$. 

Thus, the only remaining cases are $A_m$ where $m$ is odd and $11 \leq m \leq 29$. 
We observe that for $m \in \{11, 13, \dots, 29\}$ such that $m \ne 15$, we have that the only simple Weyl group $W$ of rank less than $m$ such that $\m_m (W) \ne 0$ is $A_{m - 1}$. 
Moreover, $\m_{m + 1}(W(A_{m - 1})) = 0$, so the cases of the groups $A_m$ for $m \in \{11, 13, \dots, 29\}$, $m \ne 15$ is done by considering $\m_m, \m_{m + 1}$ and $\m_{m, m+ 1}$. 

Thus, the only remaining case is that of $A_{15}$ which can be done by considering $\m_{13}, \m_{16}$ and $\m_{13, 16}$. 

This proves the theorem.
\end{proof}

We now prove Theorem \ref{thm:main2}, the main theorem of this paper. 

\begin{proof}
We recall that $W_1$ and $W_2$ denote the Weyl groups of $H_1$ and $H_2$, respectively. 
Let $m_0$ be the maximum possible among the ranks of simple factors of the groups $H_i$. 
It is clear from the \ref{cor} that $\m_i(W_1) = \m_i(W_2)$, $\m_i'(W_1) = \m_i'(W_2)$ and $\m_{i, j}(W_1) = \m_{i, j}(W_2)$ for any $i, j$. 
Then we apply Theorem \ref{thm:ind} to conclude that the product of rank $m_0$ simple factors in $W_i$ is isomorphic for $i = 1, 2$. 

Let $m$ be a positive integer less than $m_0$. 
For $i=1, 2$, let $W_i'$ be the subgroup of $W_i$ which is the product of the Weyl groups of simple factors of $H_i$ of rank $>m$. 
We assume that the groups $W_1'$ and $W_2'$ are isomorphic and then we prove that the product of the Weyl groups of rank $m$ simple factors of $H_i$ are isomorphic for $i = 1, 2$. 
This will complete the proof of the theorem by induction argument. 

Let $U_i$ be the subgroup of $W_i$ such that $W_i = U_i \times W_i'$. 
Then, since $\m_j(W_1') = \m_j(W_2')$ and $\m_j'(W_1') = \m_j'(W_2')$, we have 
$$\m_j(U_1) = \m_j(W_1) - \m_j(W_1') = \m_j(W_2) - \m_j(W_2') = \m_j(U_2) ,$$ 
$$\m_j'(U_1) = \m_j'(W_1) - \m_j'(W_1') = \m_j'(W_2) - \m_j'(W_2') = \m_j'(U_2)$$
and similarly $$\m_{i, j}(U_1) = \m_{i, j}(U_2) .$$ 

Now we use Theorem \ref{thm:ind} to conclude that the subgroups of $W_i$ which are products of the Weyl groups of simple factors of $H_i$ of rank $m$ are isomorphic for $i = 1, 2$. 

The proof of the theorem can now be completed by the downward induction on $m$. 

It also follows from the proof of the Theorem \ref{thm:ind}, that the Weyl groups of simple factors of $H_i$ are pairwise isomorphic. 
\end{proof}

\begin{rem}
We remark here that the above proof is valid even if we assume that the Weyl groups $W(H_1)$ and $W(H_2)$ share the same set of elements up to conjugacy in $GL_n(\Q)$, not just in $GL_n(\Z)$. 
Thus the Theorem $\ref{thm:main}$ is true under the weaker assumption that the groups $H_1$ and $H_2$ share the same set of maximal $k$-tori up to $k$-isogeny, not just up to $k$-isomorphism. 

We also remark that the above proof holds over the fields $k$ which admit arbitrary cyclic extensions and which have cohomological dimension $\leq 1$.
\end{rem}

\begin{rem}
Philippe Gille has recently proved $($\cite{G}$)$ that the map $\psi$, described in Lemma $\ref{lem:tori:iso}$, is surjective for any quasisplit semisimple group $H$. 
Therefore our main theorem, Theorem $\ref{thm:main}$, now holds for all fields $k$ which admit cyclic extensions of arbitrary degree. 
\end{rem}

The author thanks Joost van Hamel for informing about Gille's result.

\end{document}